\documentclass{article}
\usepackage{graphicx} 
\usepackage{amssymb,amsmath,amsthm}
\usepackage[hyphens]{url}
\usepackage{hyperref}
\usepackage{dsfont}
\usepackage{nicefrac}
\usepackage{mathtools}
\usepackage{mathrsfs}

\newcommand{\rev}{{\rm rev}}

\newcommand{\SL}{{\rm SL}}
\newcommand{\CM}{\mathcal{M}}
\newcommand{\CA}{\mathcal{A}}
\newcommand{\BZ}{\mathbb{Z}}
\newcommand{\SR}{\mathscr{R}}

\newcommand{\Prob}{\mathbb{P}}

\theoremstyle{plain}

\theoremstyle{plain}

\newtheorem{lemma}{Lemma}[section]

\newtheorem{proposition}[lemma]{Proposition}

\newtheorem{assumption}{Assumption}[section]

\theoremstyle{remark}

\newtheorem*{rem*}{Remark}
\newtheorem*{notat*}{Notation}
\newtheorem*{exm*}{Example}

\theoremstyle{definition}

\usepackage{graphicx,tikz}
\usepackage{color}

\title{Bounding the security of the Bromberg-Shpilrain-Vdovina hash function}
\author{Eilidh McKemmie}
\date{November 2026}

\begin{document}
	
	\maketitle
	
	\begin{abstract}
		The Bromberg-Shpilrain-Vdovina hash functions are a family of Cayley hash functions indexed by an integer parameter $k\ge 2$. Hashes are matrices in $SL_2(p)$ where $p$ is a large prime. We present an attack which finds collisions of length $O(\log p)$ whose average-case time complexity is $O(p^{0.28})$ when $k=2$ and $O(p^{0.36})$ when $k=3$ under an assumption in the spirit of Zaremba's Conjecture on continued fractions. This is an improvement on the more general attacks of Mullan-Tsaban and Slaughter which have average-case time complexity $O(\sqrt{p})$. For any $k$, the average-case time complexity of our attack is $O(p^{1-\delta_k+o(1)})$ where $\delta_k$ is the Hausdorff dimension of a certain Cantor set. This tends to $O(\sqrt{p})$ as $k\to\infty$, so our attack has time complexity comparable to Mullan-Tsaban and Slaughter's attacks as $k$ grows.  We give examples of collisions for large primes found by our attack when $k=2, 3$.
	\end{abstract}
	
	\section{Introduction}
	
	Let $k\ge 1$ be an integer and let
	\[A_k=\begin{pmatrix}
		1&k\\0&1
	\end{pmatrix}, \qquad B_k=\begin{pmatrix}
		1&0\\k&1
	\end{pmatrix}.\]
	
	Let $H_k:\{0,1\}^* \to \SL_2(p)$ be defined by the following:
	
	\[H(\emptyset)=I, \quad H_k(w0)=H_k(w)A_k, \quad H_k(w1)=H_k(w)B_k.\]
	
	This is a Cayley hash function. When $k=1$, this hash function was proposed by Z{\'e}mor \cite{Zemor}, and an efficient collision attack was found shortly afterwards by Tillich and Z{\'e}mor \cite{Tillich1994}. Bromberg, Shpilrain and Vdovina \cite{Bromberg_2015} proposed using the hash function with $k=2$ or $k=3$, or even larger values of $k$. This proposal has not been broken.
	
	There have been several attacks on Cayley hash functions which find collisions using palindromes or similar word reversal. Grassl et al \cite{Grassl_2010} gave the original palindrome attack on the Tillich-Z{\'e}mor hash function proposed in \cite{Tillich} (see \cite{coz2022} for another description of this attack). This attack was extended to pre-image and second pre-image attacks by Petit and Quisquater \cite{Petit_2011}. Mullan and Tsaban \cite{Mullan_2015} use word reversal to construct collisions for Cayley hash functions in $SL_2(q)$ of length $O(\log q)$ is $O(\sqrt{q})$ where $q$ is a power of $2$. Recently, Slaughter \cite{Slaughter2026} extended this result to odd prime $q$, and showed that no Cayley hash function on a pair of $2\times 2$ generators over the field of order $p$ has collision resistance greater than $O(\sqrt{p})$.
	
	Mullan-Tsaban and Slaughter's attacks \cite{Mullan_2015, Slaughter2026} find collisions of length $O(\log p)$ in $O(\sqrt{p})$ time. We present an attack which finds collisions of length $O(\log p)$ in $O(p^{1-\delta_k+o(1)})$ time assuming a Zaremba-type heuristic, where $\delta_k$ is the Hausdorff dimension of a certain Cantor set. When $k=2$ our attack has average-case time complexity $O(p^{0.28})$, and when $k=3$ it is $O(p^{0.36})$. The time complexity of our attack tends to $O(\sqrt{p})$ as $k\to\infty$, so its advantages over Mullan-Tsaban and Slaughter vanish as $k$ grows. We give examples of collisions for large primes found by our attack when $k=2, 3$. Our attack uses word reversal and, like \cite{Tillich1994}, it is a lifting attack using the fact that the Euclidean algorithm can factor target matrices which are related to a Cantor set whose dimension is known for $k=2, 3$.
	
	\section{The monoid generated by $A_k, B_k$}
	
	Let $\CM_k$ be the monoid generated by $A_k, B_k$.
	
	\begin{proposition}\label{Sanov}
		\[\CM_k \subseteq \left\{\begin{pmatrix}
			k^2a+1 & kb\\kc & k^2d+1
		\end{pmatrix} : a, b, c, d \ge 0\right\}.\]
	\end{proposition}
	\begin{proof}
		This can be proved by a simple induction argument for any value of $k$. We note that when $k=2$ the submonoid is free by Sanov \cite{sanov1947property}.
	\end{proof}

	Taking $k=1$ in the hash function construction gives a hash function proposed by Z{\'e}mor \cite{Zemor}. The monoid $\CM_1$ is the monoid of all invertible $2\times 2$ matrices with non-negative integer matrices, which leaves Z{\'e}mor's hash function open to a lifting attack which efficiently finds collisions \cite{Tillich1994}. The attack uses the fact that matrices in $\CM_1$ can be factored as products of $A_1, B_1$ using the Euclidean algorithm. This reduces the problem of finding pre-images to finding an appropriate lift of the matrix in $\SL_2(p)$ to $\SL_2(\BZ)$ which is then factored using the Euclidean algorithm. This collision attack can be turned into a pre-image attack when the attacker has an oracle for solving the integer factorization problem \cite{McKemmie2025}
	
	Bromberg, Shpilrain and Vdovina propose using $k>1$ for hashing precisely because $\CM_k$ does not have such a clean description for $k>1$. Indeed, their Problem~1 asks for a description of the shape of matrix entries in this monoid.
	
	One description of the matrices in $\CM_k$ which may prove useful involves the Euclidean algorithm. If $\begin{pmatrix}x\\y\end{pmatrix}$ appears as a column of a matrix in $\CM_k$ then it arises by applying a word in $A_k, B_k$ to either $\begin{pmatrix}1\\0\end{pmatrix}$ or $\begin{pmatrix}0\\1\end{pmatrix}$. Since \[A_k\begin{pmatrix}x\\y\end{pmatrix} = \begin{pmatrix}
		x+ky\\y
	\end{pmatrix}, \qquad B_k \begin{pmatrix}x\\y\end{pmatrix} = \begin{pmatrix}
		x\\y+kx
	\end{pmatrix},\] we observe that $\begin{pmatrix}x\\y\end{pmatrix}$ appears as a column of a matrix in $\CM_k$ exactly when the quotients produced by the Euclidean algorithm applied to $x$ and $y$ are all multiples of $k$.
	
	We can translate this into the language of continued fractions. Any rational number $\frac{x}{y}\in (0,1)$ may be written as a unique simple continued fraction $[a_1, a_2, ..., a_n]$ where the $a_i$ are positive integers and $a_n\ge 2$, such that \[[a_1, a_2, ..., a_n]=\frac{x}{y} = \cfrac{1}{a_1 + \cfrac{1}{a_2 + \cfrac{1}{a_3 + \cdots}}}.\]
	
	The $a_i$ are called the partial quotients, and are exactly the same as the quotients produced by the Euclidean algorithm applied to $x$ and $y$. Therefore, $\begin{pmatrix}x\\y\end{pmatrix}$ appears as a column of a matrix in $\CM_k$ exactly when the partial quotients in the continued fraction expansion of $\frac{x}{y}$ are all multiples of $k$.
	
	We note that this description of $\CM_k$ can also be seen by writing $A_k$ and $B_k$ as M{\"o}bius transformations \[A_k \mapsto a(z)=z+k, \qquad B_k \mapsto b(z)=\frac{z}{kz+1}.\]
	
	
	Let $E_k=\{[ka_1, ka_2, ...]\}$ be the set of all real numbers whose partial quotients are multiples of $k$. This is a Cantor set whose Hausdorff dimension $\delta_k$ is important to us, so we collect results on the dimension here.
	
	\begin{lemma}\label{lem: dim}
		With $E_k$ and $\delta_k$ as above,
		\begin{itemize}
			\item $0.71936 \le \delta_2\le 0.71950$
			\item $0.63956\le \delta_3 \le 0.64073$
			\item $\delta_k> \frac{1}{2}$ for $k\ge 2$
			\item $\delta_k\to \frac{1}{2}$ as $k\to \infty$
		\end{itemize}
	\end{lemma}
	\begin{proof}
		The first two results come from \cite[Table~1]{Chousionis_2019}.
		
		In the language of \cite[\S 2]{Chousionis_2019}, the set $E_k$ is an infinite conformal iterated function system (IFS) when $k\ge 2$, and therefore its Hausdorff dimension $\delta_k$ is the unique zero of the \textit{topological pressure} function $P$. The bounded distortion property gives $P(t)\ge \log Z_1(t)-t\log K$ where $K$ is the \textit{distortion constant} (see \cite[Proposition~2.5]{Chousionis_2019}) and $Z_1(t)=\sum_{m\ge 1}(km)^{-2t}$ so as $t\to\frac{1}{2}$ from above, $Z_1(t)\to \infty$ and therefore $P(t)\to \infty$. By \cite[Proposition~2.2]{Chousionis_2019}, $P(t)=\infty$ for $t\le \frac{1}{2}$ and $P$ is decreasing on $[0,\infty)$ and continuous on $(\frac{1}{2},\infty)$ with $P(t)\to -\infty$ as $t\to \infty$. Therefore Bowen's parameter $h:=\inf\{t\ge 0 : P(t)\le 0\}>\frac{1}{2}$, and \cite[Theorem~2.3]{Chousionis_2019} gives us $\delta_k=h>\frac{1}{2}$.
		
		Good \cite{Good_1941} proved that the Hausdorff dimension of the set of reals in $[0,1]$ whose continued fraction partial quotients are all at least $k$ tends to $\frac{1}{2}$ as $k\to \infty$. Our $E_k$ is contained in this set, so $\limsup_k\{\delta_k\}\le \frac{1}{2}$, finishing the proof.
	\end{proof}
	
	As in \cite{Tillich1994}, we will use the Euclidean algorithm to factor matrices in $\CM_k$.
	
	\begin{proposition}\label{prop: Euclid}
		If $x<y$ such that
		\begin{itemize}
			\item $\frac{x}{y}=[ka_1, ..., ka_n]\in E_k$
			\item $y \equiv 1 \text{ (mod }k^2)$
			\item $x \equiv 0 \text{ (mod }k)$
		\end{itemize} then \[B_k^{a_1}A_k^{a_2}B_k^{a_3}\cdots = \begin{pmatrix}
			*&x\\ *&y
		\end{pmatrix}.\]
	\end{proposition}
	\begin{proof}
		The congruence conditions come from Proposition~\ref{Sanov}, and the matrix identity follows from the definitions of $A_k$ and $B_k$.
	\end{proof}
	
	\section{A collision attack}
	\subsection{The reverse of a word}\label{sec: reverse}
	
	For $w\in \{0, 1\}^*$, define $\rev(w)$ to be the reverse of $w$. If $C=\begin{pmatrix}
		0&1\\1&0
	\end{pmatrix}$ then $A=CA^TC$ and $B=CB^TC$. Therefore, $H_k(\rev(w))=CH_k(w)^TC$.
	
	If we know a word $w$ such that $H_k(w)=\begin{pmatrix}
		a&b\\c&dp
	\end{pmatrix}$, then \begin{align*}H_k(w)H_k(\rev(w))&= \begin{pmatrix}
			a&b\\c&dp
		\end{pmatrix}\begin{pmatrix}
			dp&b\\c&a
		\end{pmatrix} \\&=
		\begin{pmatrix}
			adp+bc&2ab\\2cdp&adp+bc
		\end{pmatrix}\equiv \begin{pmatrix}
			-1&2ab\\0&-1
		\end{pmatrix} \, \text{ mod }p.\end{align*}
	
	This matrix commutes with $A_k$, which gives the following collision \[H_k(w\rev(w)0)= H_k(0w\rev(w)).\]
	
	Note that $w\rev(w)0\ne 0w\rev(w)$ unless every letter of $w$ is a $0$.
	
	\subsection{Finding appropriate lifts}
	
	It remains to find a reliable algorithm for finding a word $w$ such that $H_k(w)=\begin{pmatrix}
		a&b\\c&dp
	\end{pmatrix}\in \CM_k$. By Proposition~\ref{Sanov}, we must have $dp \equiv 1 \text{ (mod }k^2)$, so we fix a small value of $d$ such that $dp \equiv 1 \text{ (mod }k^2)$. The size of this value of $d$ will be discussed in Section~\ref{sec: d-switch}.
	
	Start with a prefix word $w_1$ so that $H_k(w_1)\in \CM_k$ by construction, and denote by $(c_1, d_1)$ the bottom row of $H_k(w_1)$. We pick $w_1$ at random such that $c_1, d_1 \sim \sqrt{p}$.
	
	If $w_1$ has length $n$ then $c_1, d_1 \sim e^{n\lambda_k}$ where $\lambda_k$ is the growth rate of a random product of $\{A_k, B_k\}$. This $\lambda_k$ is called the Lyapunov exponent, and is of independent interest \cite{Pollicott_2010}.
	
	\begin{lemma}
		If $w_1$ is such that $c_1, d_1 \sim \sqrt{p}$ then the expected length of $w_1$ is $O(\log p)$.
	\end{lemma}
	\begin{proof}
		The group $\langle A_k, B_k \rangle$ is Zariski dense in $SL_2(\mathbb{R})$ since its Zariski closure contains all upper and lower unitriangular matrices, and must therefore be $SL_2(\mathbb{R})$. The measure $\frac{1}{2}\left(\delta_{A_k}+\delta_{B_k}\right)$ has finite support, so we may use Furstenberg's Law of Large Numbers as stated in \cite[Theorem~0.6]{Benoist_2016}, which shows that $\lambda_k$ exists \cite{Furstenberg_1960} and $\lambda_k>0$ \cite[Theorem~8.6]{Furstenberg_1963}. Therefore for a fixed $k$, we expect $w_1$ to have length $\sim \frac{\log{p}}{2\lambda_k}=O(\log p)$.
	\end{proof}
	
	Now we consider solutions $x, y \sim \sqrt{p}$ to the linear equation $c_1x+d_1y=dp$. If $\frac{x}{y}\in E_k$, then the Euclidean algorithm gives us a word $w_2$ such that $H_k(w_2)$ has last column $\begin{pmatrix}
		x\\y
	\end{pmatrix}$, as described in Proposition~\ref{prop: Euclid}. Now the matrix $H_k(w_1w_2)$ has $(2, 2)$-entry $c_1x+d_1y=dp$, which gives us a collision as described in \S\ref{sec: reverse}. If $\frac{x}{y}\notin E_k$ then we try again with another random prefix $w_1$.
	
	The remainder of this section discusses how many random prefixes $w_1$ we must try on average.
	
	Let $S(w_1)$ be the number of times $w_1$ appears as a successful prefix
	\begin{align*}S(w_1)&:= \#\left\{(x, y) : \frac{x}{y}\in E_k, x, y \le \sqrt{p}, c_1x+d_1y=dp\right\}.\end{align*}
	
	We continue choosing random prefixes $w_1$ until we find one with $S(w_1)\ge 1$. Since this follows the geometric distribution, the expected number of prefixes we must check is $1/\Prob(S\ge 1)$.
	
	We will make the following assumption.
	
	\begin{assumption}\label{assumption}
		Assume as $w_1$ ranges over words of length $O(\log p)$ that $\frac{x}{y}$ is in $E_k$ as frequently as a random rational $\frac{a}{b}$ with $a, b \le \sqrt{p}$.
	\end{assumption}
	
	Zaremba's conjecture, which states that every positive integer denominator has at least one coprime numerator so that the partial quotients are all bounded by a small constant, has been proved for a set of denominators of density $1$ \cite{Bourgain_2014, kontorovich2013apollonius} and recently announced in full generality by Zhang \cite{Zhang2026} building upon work of Shkredov \cite{Shkredov2026}. A generalization called the Local-Global Conjecture \cite{Bourgain_2014} has been proven false in several cases \cite{Rickards_2025, Shi2025}.
	
	The assumption we are making is purely conjectural, and is in the spirit of the Zaremba Conjecture but for an infinite alphabet. In Section~\ref{sec: experiments}, we explain what to do when finding appropriate $x, y$ is more difficult than expected, and provide experimental evidence that the algorithm works well.
	
	\begin{proposition}\label{prop: count}
		Under Assumption~\ref{assumption}, for any $\epsilon>0$ we have $\Prob(S\ge 1)\gg p^{\delta_k-1-\epsilon}$.
	\end{proposition}
	
	First, we write down the density of $E_k$ in $[0, 1]$ using the notation of the survey \cite{kontorovich2013apollonius}. Let
	\[\SR_{\CA}(N):=\left\{\frac{x}{y} = [a_1, a_2, ..., a_n] : (x, y)=1, \  1\le x < y < N, \ a_i\in \CA, \ \forall i \right\}.\]
	
	\begin{lemma}\label{lem: density}
		For every $\epsilon >0$, \[\#\SR_{E_k}(N) \gg N^{2\delta_k - \epsilon}.\]
	\end{lemma}
	\begin{proof}
		Hensley \cite{hensley1989distribution} (see also \cite[(2.17)]{kontorovich2013apollonius}) showed that, when $\CA$ is a finite set and $\delta_\CA$ is the Hausdorff dimension of $\{[a_1, a_2, ...]: a_i\in \CA, \forall i\}$, we have $\#\SR_{\CA}(N) \asymp N^{2\delta_\CA}$. By \cite[Theorem~2.3]{Chousionis_2019}, $\delta_k = \sup\{\delta_\CA : \CA \subset E_k \text{ finite}\}$ so we get $\#\SR_{E_k}(N) \gg N^{2\delta_k - \epsilon}$.
	\end{proof}
	
	\begin{proof}[Proof of Proposition~\ref{prop: count}]
		There are around $\frac{d}{k}=O(1)$ positive integer points on the line $c_1x+d_1y=dp$, and under Assumption~\ref{assumption} we may invoke Lemma~\ref{lem: density} to see that these each have probability $\gg \sqrt{p}^{(2\delta_k-2\epsilon)}/\sqrt{p}^2$ of having $\frac{x}{y}\in E_k$.
	\end{proof}
	
	Therefore for any $\epsilon>0$ the expected number of prefixes $w_1$ we must check is at most $1/\Prob(S\ge 1)=O(p^{1-\delta_k+\epsilon})$, the average-case time complexity of our attack under Assumption~\ref{assumption}.

	\section{Implementation and experimental results}\label{sec: experiments}
	
	We used Claude Opus 4.8 to implement the attack. Since our algorithm relies on Assumption~\ref{assumption}, we build in a few safeguards to ensure our word length and run time do not become too large. We describe these safeguards first, then include the results of testing which shows the algorithm performs with the claimed time complexity. Finally, we give examples of collisions for large primes which have not been found before.
	
	\subsection{Capping the length}
	
	If Assumption~\ref{assumption} is true then the algorithm as described gives words of length $O(\log p)$ in time $O(p^{1-\delta_k+o(1)})$, but if the assumption fails the algorithm may return longer words. To ensure the algorithm will only return a collision when it is short enough, we add a step which rejects any $w_2$ with $|w_2| > 8\log_2 p$. Note that this restricts us to a finite subset of $E_k$ which has dimension strictly smaller than $\delta_k$. This cap on the length of $w_2$ was chosen after computational experiments, and following sections show that the algorithm works well using this cap.
	
	\subsection{The size of $d$}\label{sec: d-switch}
	Whether or not Assumption~\ref{assumption} is close to the truth depends on our choice of $d$, so we change $d$ if no appropriate lift is found after the expected number of trials.
	
	The first choice of $d$ is $0< d \le k^2$, and $d$ is increased by $k^2$ if no appropriate lift is found after $2^{\lfloor (1-\delta_k)n\rfloor}$ words $w_1$ have been tried, where $n$ is the number of bits of the prime $p$.
	
	\subsection{The algorithm}
	\begin{enumerate}
		\item Choose $0< d \le k^2$ such that $dp \equiv 1 \text{ (mod }k^2)$.
		\item Choose a random word $w_1$ so that if $(c_1, d_1)$ is the bottom row of $H_k(w_1)$ then $c_1, d_1\sim \sqrt{p}$. The word $w_1$ has length $O(\log p)$.
		\item Solve the linear equation $c_1x+d_1y=dp$ in the integers.
		\item Apply the Euclidean algorithm to determine whether $\frac{x}{y}\in E_k$.
		\item If $\frac{x}{y}\notin E_k$, return to Step 2. Otherwise, construct $w_2$ such that $H_k(w_2)=\begin{pmatrix}
			*&x\\ *&y
		\end{pmatrix}$ using Proposition~\ref{prop: Euclid}.
		\item If $|w_2| > 8\log_2 p$, return to Step 2. Otherwise, we have found $w:= w_1w_2$ such that $H_k(w \rev(w)0)=H_k(0w \rev(w))$, a collision.
		\item If we try more than $2^{\lfloor (1-\delta_k)n\rfloor}$ choices of $w_1$ without success, increase the value of $d$ by $k^2$ and return to Step 2.
	\end{enumerate}
	
	Under Assumption~\ref{assumption}, this algorithm returns a collision of length $O(\log p)$ in $O(p^{1-\delta_k+o(1)})$ time, on average. By Lemma~\ref{lem: dim}, this gives average-case time complexity of $O(p^{0.28})$ when $k=2$, and $O(p^{0.36})$ when $k=3$. As $k$ grows, the average-case time complexity tends to $O(\sqrt{p})$, and so the advantages of our attack over Slaughter's attack \cite{Slaughter2026} vanish.
	
	\subsection{Experimental evidence for time complexity}
	We implemented the algorithm in Python and ran it on 60 random primes per bit-length for $k=2, 3$. The tables below record the mean and median number of prefixes $w_1$ tried until the algorithm returned a collision.
	
	\subsubsection{$k=2$}
	For $k=2$, the least-squares exponent of the mean trial count is $0.297$, close to the
	predicted $1-\delta_2 \approx 0.28$.
	\medskip
	\begin{center}
		\begin{tabular}{r|r|r}
			$\log_2 p$ & mean trials & median trials \\
			\hline
			20 &      17.9 &     13  \\
			25 &      57.1 &     33  \\
			30 &     171.3 &    144  \\
			35 &     334.3 &    262  \\
			40 &     793.6 &    576 \\
			45 &    3502.2 &   2429  \\
			50 &    8862.5 &   5028  \\
			55 &   27613.1 &  19240   \\
			60 &   68548.5 &  36156  \\
		\end{tabular}
	\end{center}
	
	\subsubsection{$k=3$}
	For $k=3$, the least-squares exponent of the mean trial count is $0.365$, close to the
	predicted $1-\delta_3 \approx 0.36$.
	\medskip
	\begin{center}
		\begin{tabular}{r|r|r}
			$\log_2 p$ & mean trials & median trials  \\
			\hline
			20 &      56.8 &     32 \\
			25 &     221.9 &     94   \\
			30 &     651.6 &    493  \\
			35 &    2248.4 &   1200  \\
			40 &    6408.8 &   5302 \\
			45 &   29010.3 &  23246  \\
			50 &  136093.4 & 101010 \\
		\end{tabular}
	\end{center}
	
	\subsection{Examples of collisions for large primes}
	To find collisions for large primes, we implemented the algorithm in C using the GNU Multiple Precision Arithmetic Library (GMP). This runs much more efficiently than our Python implementation.
	
	Using a 2-vCPU cloud container (Intel Xeon, ~2.3 GHz, 8 GB RAM), we found collisions for $k=2$ with a $112$-bit prime, and for $k=3$ with a $92$-bit prime. For $k=2$ with a $152$-bit prime, we found a collision using a Google Compute Engine with 88 vCPUs (Intel Xeon Platinum 8481C, ~2.7 GHz). For $k=3$ with a $118$-bit prime, we found a collision using a Google Compute Engine with 44 vCPUs (Intel Xeon Platinum 8481C, ~2.7 GHz).
	
	All collisions were independently checked using a simple Python script. The colliding words are too long to print here but are included in the Python script available with the arXiv version of this paper.
	
	As far as the author knows, these are the largest primes for which collisions for these hash functions have been found. The number of words $w_1$ tried in each case is listed in the following table as well as the length of the collision (i.e. the sum of the lengths of the colliding words, at most $2+32\log_2 p$).
	
	\medskip
	\begin{center}
		\begin{tabular}{r|r|r|r|r|r}
			$k$& $p$ & length & trials & $\approx p^{1-\delta_k}$ & vCPUs\\
			\hline
			\rule{0pt}{3ex}2 & $2^{112}-75$ & 878 & $2.9\times 10^9$ & $2.7 \times 10^9$&2\\
			2 & $2^{152}-17$ & 3626 & $2.3\times 10^{10}$ & $6.5\times 10^{12}$ & 88\\
			\hline
			\rule{0pt}{3ex}3 & $2^{92}-83$ & 642 & $1.7 \times 10^8$ & $4.9\times 10^9$ & 2\\
			3 & $2^{118}-5$ & 2494 & $6.1\times 10^{10}$ & $2.7\times 10^{12}$ & 44\\
		\end{tabular}
	\end{center}

	\subsection*{Code availability}
	The code for this project is available in an ancillary file in the arXiv version of this paper. It contains the Python implementation of the algorithm 
	\texttt{bsv\_search.py}, the faster version written in C \texttt{search.c}, and the Python code verifying the reported large-prime collisions \texttt{large\_p\_verification.py}.
	
	\subsection*{Acknowledgements}
	The author thanks Chris Lutsko for interesting discussions, and in particular pointing out that the quotients of columns of positive products of $A_2$ and $B_2$ are those whose continued fraction expansions have even partial quotients. The author thanks Vladimir Shpilrain for comments on an early draft of this paper.
	
	\subsection*{AI use acknowledgement}
	Aristotle by Harmonic was used in exploratory conversations and suggested the main approach of this paper, along with others that did not prove fruitful. The author found a significant gap in the initial suggestion and prompted Aristotle to patch the gap using Assumption~\ref{assumption}. Claude Opus 4.8 wrote the code which produced the experimental results, and was used to teach the author about Zaremba's conjecture, proofread and check the article before submission for peer review. No LLM-generated writing appears in this paper. The author takes responsibility for the entire contents of this paper.
	
	\bibliographystyle{plain}
	\bibliography{bib}
	
\end{document}